\documentclass[12pt]{amsart}
\usepackage{amsfonts,graphics,amsmath,amsthm,
amsfonts,amscd, amssymb,amsmath,latexsym}
\usepackage{epsfig}

\theoremstyle{plain}
\newtheorem*{corollary*}{Corollary}

\newtheorem{theorem}{Theorem}
\newtheorem*{theorem*}{Theorem}
\newtheorem*{lemma*}{Lemma}
\newtheorem{lemma}{Lemma}

\newtheorem*{claim*}{Claim}

\newtheorem*{example*}{Example}

\newtheorem*{proof*}{proof}
\newtheorem*{prop*}{Proposition}
\newtheorem{prop}{Proposition}
\newtheorem{Definition/Proposition}{Definition/Proposition}
\def\NE{\overline{NE}}
\newtheorem*{question*}{Question}

\newtheorem{definition}{Definition}
\def\M{{\overline{\mathnormal{M}}}}
\def\S{\mathnormal{S}}

\def\T{\mathnormal{T}}

\begin{document}
\title{\textbf{The Mori cones of moduli spaces of pointed curves \\
of small genus }}
\author{ GAVRIL FARKAS AND ANGELA GIBNEY}
\maketitle

\vspace{2ex} \normalsize
\section{\textbf{Introduction}} 

In this paper we describe the Mori cone of curves of the moduli space
 $\M_{g,n}$ of $n$-pointed stable curves of small genus $g$. Although 
 important aspects of the birational geometry of $\M_{g,n}$ heavily depend 
 on whether $g$ is large with respect to $n$ (e.g. the Kodaira dimension), 
 it turns out that surprisingly the problem of determining the cone of curves 
 $\NE_1(\M_{g,n})$ can be expressed entirely in terms of the geometry of the moduli 
 spaces  $\M_{0,g+n}$ of rational curves with marked points (cf. \cite{GKM}).

There is a combinatorial stratification of $\M_{g,n}$ given by topological type and the
components of the $1$-dimensional stratum, that is, loci of curves with $(3g-4+n)$-nodes are
called {\sl Faber curves} (or $F$-curves). Our first result is that any curve in $\M_{g,n}$ is
numerically equivalent to an effective combination of these $F$-curves at least when the genus is relatively
small:

\begin{theorem}\label{ne1}
The Mori cone of curves $\NE_1(\M_{g,n})$ is generated by $F$-curves when $g\leq 13,n=0$
or $g\leq 8,n=1$ or $g=6,n=2$.
\end{theorem}

Thus in this range the cone of nef divisors is described by a simple set of inequalities corresponding
to the numerical properties of all $F$-curves (cf. \cite{GKM}). This result was
known when $n=0$ for $g\leq 11$ (cf. \cite{KMcK}), and when $n=1,g\leq 6$ (cf.
\cite{GKM}). The case $n=0,g\leq 4$ was first settled by Faber (cf. \cite{Fa1}).

Already for small $n$ the moduli spaces $\M_{0,n}$ are quite intricate objects  deeply rooted in
classical algebraic geometry. For instance $\M_{0,5}$ is a del Pezzo surface of degree $5$ while
$\M_{0,6}$ is a small resolution of two famous modular varieties: the {\sl Segre cubic}
$\mathcal{S}_3\subseteq \mathbb P^4$ which is the symmetric GIT moduli space of $6$ points on
$\mathbb P^1$ and is the unique cubic threefold with $10$ nodes, and its dual, the {\sl Igusa
quartic} $\mathcal{I}_4\subseteq \mathbb P^4$ which is the GIT moduli space of $6$ points on
$\mathbb P^2$ which lie on a conic and whose singular locus consists of $15$ double lines (cf.
\cite{H}).

Since the space $\M_{0,n}$ has a combinatorial description somewhat similar to that of a 
toric variety (although $\M_{0,n}$ itself is certainly not a toric variety), Fulton asked
whether any effective (nef) divisor on $\M_{0,n}$ is linearly equivalent to an effective 
combination of boundary divisors corresponding to singular curves. We prove the following result:

\begin{theorem}\label{fulton} The cone of nef divisors of $\M_{0,6}$ is contained in the convex hull of boundary classes and it has a natural decomposition into $11$ subcones.
 \end{theorem}

The precise inequalities defining these subcones can be found in Section 4. The explicit decomposition of the nef cone of $\M_{0,6}$ enables us to 
classify all fibrations of
$\M_{0,6}$. Recall that a morphism $f:X\rightarrow Y$ with $X$ and $Y$ being irreducible
projective varieties, is said to be a {\sl fibration} if $\mbox{dim}(X)>\mbox{dim}(Y)$ and $f_*\mathcal{O}_X=\mathcal{O}_Y$, that is, $f$ is its own Stein factorization. We have the following:
\begin{theorem}\label{fibrations}
1. Any fibration of $\M_{0,5}$ factors through a projection $\M_{0,5}\rightarrow \M_{0,4}$ dropping one of the marked points.

2. Any fibration of $\M_{0,6}$ factors through a projection $\M_{0,6}\rightarrow \M_{0,n}, n\in \{4,5\}$, dropping one or two points,  or through a projection $\M_{0,6}\rightarrow \M_{0,4}\times \M_{0,4}$ obtained by dropping two disjoint pairs of marked points.
\end{theorem}
The fact that every nef divisor on $\M_{0,6}$ is linearly equivalent to an effective combination of boundary divisors (that is, the first part of Theorem \ref{fulton}) has been previously checked by Faber (cf. \cite{Fa2}) and by Keel (using the computer program Porta). The salient features of our Theorem \ref{fulton} are the method of proof which can be applied in more general situations (see Propositions 8 and 9) and the decomposition of the nef cone of $\M_{0,6}$ into subcones which in particular leads to a classification of all fibrations of $\M_{0,6}$. For example the fibrations $\M_{0,6}\rightarrow \M_{0,4}\times \M_{0,4}$ correspond to nef divisors lying  in the boundary of two different chambers (see Section $4$ for a precise statement). We also mention that Theorem \ref{fulton} gives a new combinatorial proof that $\NE_1(\M_{0,6})$ is spanned by $F$-curves (cf. \cite{KMcK}, Theorem 1.2). The original proof used that $\M_{0,n}$ is a $\mathbb Q$-Fano variety for $n\leq 6$
(which is not the case for any $n\geq 7$).

We note that Sean Keel showed that there are effective divisors on $\M_{0,6}$ which are not  expressible as
effective combinations of boundary classes (see also \cite{Ve}). For example, if $\phi:\M_{0,6}\rightarrow
\M_{3}$ is the map obtained by identifying three pairs of points on a $6$-pointed rational curve and
$\overline{\mathcal{H}}\subseteq \M_3$ is  the locus of hyperelliptic curves then
$\phi^*(\overline{\mathcal{H}})$ is not linearly  equivalent to an effective sum of boundary classes.  
Hassett and Tschinkel recently proved that the effective cone on $\M_{0,6}$ is spanned by the boundary
classes and the pullbacks $\phi^*(\overline{\mathcal{H}})$  corresponding to all possibilities of
identifying three pairs of points (cf. \cite{ht}).  In light of their work, Theorem \ref{fulton} shows that the nef cone
of $\M_{0,6}$ is very small with respect to the effective cone of divisors.  Indeed, while we show that any nef divisor may
be expressed as an effective sum of the $16$ boundary classes, one needs another   $15$ divisor classes to describe all effective
divisors.

Our proof of Theorem \ref{ne1} makes use of the following bridge theorem of Gibney, Keel and
Morrison (cf. \cite{GKM}, Theorem 0.3): if $\psi:\M_{0,g+n}/S_g\rightarrow \M_{g,n}$ is the
map given by attaching elliptic tails to the first $g$ marked points of a $(g+n)$-pointed
rational curve, then a divisor $D$ on $\M_{g,n}$ is nef if and only if $\psi^*(D)$ is nef
and $D$ meets all $F$-curves on $\M_{g,n}$ nonnegatively. In other words, to show that
$\NE_1(\M_{g,n})$ is generated by $F$-curves it suffices to prove the similar statement on
the space $\M_{0,g+n}/S_g$.

We note that Theorem \ref{fibrations} should be compared to Gibney's result that for $g\geq 2$ any fibration
of $\M_{g,n}$ factors through a projection to some $\M_{g,i}$ ($i<n$) dropping some of the marked points
(cf. \cite{g},\cite{GKM}, Corollary 0.10). Paradoxically, because of the combinatorial complexity of $\mbox{Pic}(\M_{0,n})$, the fibration problem is much more difficult in genus $0$ than in higher genus!

\vskip 6pt
\textbf{Acknowledgments:} We are grateful to Igor Dolgachev and Sean Keel for many interesting discussions related to this project.

\newpage 
\section{\textbf{Generalities on $\M_{0,n}$}}

We record a few facts about the moduli space $\M_{0,n}$ of stable rational $n-$pointed
curves.  For more information about $\M_{0,n}$ see for example  ~\cite{Kap} or \cite{Ke}. Throughout the paper we work exclusively with $\mathbb Q$-divisors and all the Picard groups we consider are with rational coefficients.

A {\sl vital codimension-$k$-stratum} is a component of the closure of the locus of points in
$\M_{0,n}$ that correspond to curves with $k$ nodes.  The boundary of $\M_{0,n}$ is composed
of the vital codimension $1$-strata $\Delta_{S}$ where $S\subset \{1,\ldots,n\}$ with
$|S|,|S^{c}|\geq 2$. We denote by $\delta_{S}$ the linear equivalence class
of $\Delta_{S}$ in $\rm{Pic}$$(\M_{0,n})$.   An effective $1$-cycle that is numerically equivalent to a vital $1$-stratum is also known as an {\sl F-curve}.  By an {\sl
F-divisor} we mean a divisor than nonnegatively intersects the $F$-curves.

We will consider the tautological classes $\psi_i=c_1(\mathbb L_i)$ for $1\leq i\leq n$, where
$\mathbb L_i$ is the line bundle on $\M_{0,n}$ whose fibre over the moduli point
$[C,x_1,\ldots, x_n]$ is $T_{x_i}^{\vee}(C)$. Recall also that there exists an ample divisor class
$$\kappa_1=\sum_{\stackrel{S\subset \{1,\ldots,n\}}{|S|\leq n/2}} \frac{(|S|-1)(n-|S|-1)}{(n-1)} \delta_S$$ whose support is the whole boundary of $\M_{0,n}$ (cf. \cite{AC}).

For each subset $S\subset \{1,\ldots,n\}$ we denote by $G_S$ the stabilizer of $\delta_S$ under the natural
action of $\S_n$ on $\M_{0,n}$. Then the $G_S$-invariant divisor classes of the form
$$
\delta_{b}^{\S,a}:=
\sum_{
          \stackrel{A \subset \S , |A| = a }
         {B \subset \S^c , |B| = b}}^{}
  \delta_{A \cup B}
$$
generate $\mbox{Pic}(\M_{0,n})^{G_S}$. We have the following relation between tautological and boundary classes:

\begin{lemma}\label{psi}

The tautological classes $\psi_i$, for $1 \leq i \leq n$ have  the following average
expression in terms of $G_{\{i\}}$-invariant boundaries: \[ \psi_i \equiv \sum_{j=1}^{n-3} \space  \frac{(n-1-j)(n-2-j)}{(n-1)(n-2)} 
\delta^{\{i\},1}_j. \] 
\end{lemma}

 \begin{proof}  We use that given two distinct elements $q,r\in \{1,
\ldots, n\}-\{i\}$, we have that $\psi_i \equiv \sum_{\stackrel{i \in \S}{q,r \notin \S}}^{}\delta_{\S}$ (cf. \cite{AC}, Proposition 1.6). We then average all such relations obtained by varying $q$ and $r$.
 \end{proof}

The following average relation between $G_S$-invariant divisor classes will be used throughout the paper:

 \begin{prop}\label{average}
Suppose $\S \subset \{1,\ldots,n\}$ has $s$ elements. The following
  relation in
 $\rm{Pic}$$(\M_{0,n})$ holds:
\begin{eqnarray*}
\delta_S=\sum_{ \stackrel{1 \leq a \leq s,(a,b)\neq (s,0)}{0 \leq b \leq n-s-1}}^{} \eta_{s,a,b}
 \ \delta_{b}^{\S,a},
  \end{eqnarray*} 
where 
\begin{eqnarray*}
\eta_{s,a,b} & := &
\frac{a(b+s-n)\bigl(1+b+a(n-1)-n+s-s(a+b)\bigr)}{s(s-1)(n-s)(n-s-1)}.
\end{eqnarray*} 
\end{prop}

\begin{proof}
 We use Keel's relation in $\mbox{Pic}(\M_{0,n})$ (cf. \cite{Ke}):  given
  four distinct elements $p,q,r,s \in \{1, \ldots, n\}$ we have that
\[
\sum_{\stackrel{p,q \in \T}{r,s \notin \T}}^{} \delta_{\T} \equiv
  \sum_{\stackrel{p,r \in \T}{q,s \notin \T}}^{} \delta_{\T} \equiv
  \sum_{\stackrel{p,s \in \T}{q,r \notin \T}}^{} \delta_{\T}.
\]
Having fixed $S$ we write down all possible such relations for which $p,q\in S$ and $r,s\notin S$. Then we add them together and average.  
\end{proof}

It is well known that the boundary classes $\delta_S$ generate $\mbox{Pic}(\M_{0,n})$ (cf.
\cite{Ke}). The existence of many relations between the $\delta_S$'s, hence the absence of a
\lq \lq canonical" basis of $\mbox{Pic}(\M_{0,n})$ reflects the combinatorial
complexity of $\M_{0,n}$. Using Kapranov's description of $\M_{0,n}$ as the space obtained
from $\mathbb P^{n-3}$ after a sequence of $(n-4)$ blow-ups one sees that $\psi_n$ and the
boundaries $\delta_{S\cup\{n\}}$ with $S\subset \{1,\ldots, n-1\}$ and $|S|\leq n-4$, constitute a basis for
$\mbox{Pic}(\M_{0,n})$. In particular $\rho(\M_{0,n})=2^{n-1}-{n\choose 2}-1$. However, this
basis singles out the $n$-th marked points and we chose to express all our calculations in a
basis which treats all marked points equally: 

\begin{lemma}\label{basis} For $n \geq 5$, the classes
$\{\psi_i\}_{i=1}^n$ and $\delta_S$ where $|S|,|S^c|\geq 3$ form a basis of
$\rm{Pic}(\M_{0,n})$. 
\end{lemma} 
\begin{proof} We denote by $V\subset \mbox{Pic}(\M_{0,n})$ the subspace generated by the classes $\{\psi_i\}_{i=1}^n$ and $\{\delta_S\}_{|S|,|S^c|\geq 3}.$  It is enough to show that $\delta_{xy}\in V$ for all distinct $x,y$. From Proposition \ref{average} we obtain that $ (n-2)(n-3)\delta_{xy}+2\delta_{2}^{xy,0}-(n-3)\delta_{1}^{xy,1}\in V\mbox{ }(i)$. 

By writing the relation $\psi_i=\sum_{\stackrel{i\in T}{x,y\notin T}} \delta_T$
for all $i\in \{x,y\}^c$ and averaging we obtain that $(n-2)\delta_{xy}+2\delta_2^{xy,0} \in V \mbox{ }(ii)$. Finally, by averaging all relations 
$\psi_x+\psi_y=\sum_{\stackrel{x\in T}{a,b\notin T}}\delta_T+\sum_{\stackrel{y\in T}{a,b\notin T}}\delta_T$ over all $a,b\in \{x,y\}^c$ we obtain that $$ 2{n-2\choose 2}\delta_{xy}+2\delta^{xy,0}_2+{n-3\choose 2} \delta_{1}^{xy,1}\in V \mbox{ }\mbox{ }(iii).$$ Clearly $(i)-(iii)$ imply that $\delta_{xy}\in V$.
\end{proof}

We will often use the following notation:

\begin{definition} For a divisor $
D \equiv \sum_{1 \leq i \leq n}^{}c_i \psi_i - \sum_{|S|,|S^c|\geq 3}^{}b_{\S}\delta_{\S}$ on $\M_{0,n}$
 and for a fixed subset $\T \subset \{1,\ldots,n\}$,  we set
$$I^{\T} := \sum_{t \in \T}c_t, \ \ O^{\T} := \sum_{t \notin \T}c_t, \ \
\Sigma_{i}^{\T,j} := 
\sum_{\stackrel{A \subset \T , |A| = j}{B \subset \T^c, |B| = i}}^{}
b_{A \cup B}.$$
\end{definition}

We also recall that $F$-curves in $\M_{0,n}$ correspond to  partitions 
$I,J,K,L$ of $\{1,\ldots,n\}$ into non-empty subsets. For each such partition we
have a map $\nu:\M_{0,4}\rightarrow \M_{0,n}$ obtained by attaching $1+|I|,1+|J|,1+|K|$ and $1+|L|$-pointed rational curves at each of the
four marked points. Every $F$-curve in $\M_{0,n}$ is numerically equivalent 
to such an image $\nu(\M_{0,4})$ corresponding to a partition (cf. \cite{GKM}, Theorem 2.2).

\section{\textbf{The fibrations of $\M_{0,5}$}}

In this section we first show that any $F$-nef divisor in $\M_{0,5}$ can be expressed as an
effective sum of boundary classes.  Although this result can be  proved in various ways we
present it because it illustrates our technique for giving a natural presentation of any
divisor in  terms of boundary classes via averaging. Moreover, it enables us to classify the
fibrations of $\M_{0,5}$.

For $a,b \in \{1,\ldots,5\}$ we consider the $G_{ab}$-invariant sum of $F$-curves
$C^{ab}:=\sum_{i\in \{a,b\}^c} \Delta_{abi}$. We show that any divisor on $\M_{0,5}$ has a
canonical presentation in terms of boundary divisors.

\begin{prop}\label{m05}

If $D$ is any divisor in $\M_{0,5}$  then  $$ D \equiv \sum_{a,b \in \{1,\ldots,5\}}^{}
\frac{1}{6} \ \ \bigl(C^{ab} \cdot D  \bigr) \ \  \delta_{ab}. $$ In particular any
$F$-divisor is an effective sum of boundary classes. \end{prop}

\begin{proof}
We have seen that $\{\psi_i \}_{i=1}^5$ forms a basis for $\mbox{Pic}(\M_{0,5})$.   Let $D \equiv  \sum_{1 \leq i \leq 5}^{}c_i \psi_i $ be any divisor on $\M_{0,5}$. Using the average formula from Lemma \ref{psi}
$$\psi_i=\frac{1}{2}(\sum_{a\neq i}\delta_{ai})+\frac{1}{6}(\sum_{a,b\neq i} \delta_{ab})$$ 
we can rewrite $D$ as
$$D \equiv  \sum_{a,b \in \{1,\ldots,5\}}^{} \frac{1}{6}\bigl( 3\  I^{ab} + O^{ab} \bigr) \delta_{ab}.
$$
The coefficient of $\delta_{ab}$  
is this expression is just $\frac{1}{6} ( D \cdot C^{ab} )$ so the conclusion follows.
\end{proof}
\noindent \textbf{Remark.} If $D\equiv \sum_{i=1}^5 c_i\psi_i$ is an $F$-divisor on $\M_{0,5}$ we see that $D\cdot \Delta_{ab}=O^{ab}\geq 0$ for any $a,b\in \{1,\ldots,5\}$. Moreover if $D\cdot C^{ab}=3I^{ab}+O^{ab}=0$ then $c_i=-(c_a+c_b)\geq 0$, for all $i\in \{a,b\}^c$.

Next we prove that a nontrivial 
$F$-divisor on $\M_{0,5}$ is either big or the pull-back of an ample divisor under the projection $\pi_i:\M_{0,5}\rightarrow \M_{0,4}$ dropping the $i$-th point.

\vskip 4pt
\noindent
\emph{Proof of Theorem \ref{fibrations}, Part 1}. Let $D\equiv \sum_{i=1}^5 c_i\psi_i$ be
a nontrivial $F$-divisor. We have the following possibilities: \newline \indent  \textbf{1.}
$D\cdot C^{ab}>0$ for any $a,b\in \{1,\ldots,5\}$. Then using Proposition \ref{m05} we can
write $D\equiv a \kappa_1+(\mbox{Effective})$, for some $a \in \mathbb Q_{>0}$
and since $\kappa_1$ is ample $D$ has to be big so it does not give rise to a fibration.
Thus we may assume that $D\cdot \delta_{ab}=0$ for some $a,b\in \{1,\ldots,5\}$, say $D\cdot
C^{12}=0$ . There are two possibilities:

\indent \textbf{2.} $D\cdot C^{1i}>0$ for all $i\in \{1,2\}^c$. Then
$c_3=c_4=c_5=c=-(c_1+c_2)>0$. Moreover $D\cdot C^{ab}>0$ for $a,b\in \{1,2\}^c$
and $D\cdot C^{2i}>0$ for all $i\in \{1,2\}^c$. In this
case using Proposition \ref{m05} the divisor $D$ can be rewritten as a {\sl positive} combination  $$D\equiv
\frac{5c+3c_1+c_2}{6}D_1+\frac{5c+c_1+3c_2}{6} D_2,$$ where $D_1=\sum_{a,b\in
\{1,2\}^c}\delta_{ab}+\sum_{a\neq 1,2}\delta_{1a}$ and $D_2=\sum_{a,b\in \{1,2\}^c}
\delta_{ab}+\sum_{a\neq 1,2} \delta_{2a}.$

Since $D_2=\pi_1^*(\delta_{23}+\delta_{24}+\delta_{34})$ and
$D_1=\pi_2^*(\delta_{13}+\delta_{14}+\delta_{34})$, it follows that $D$ is the pull-back of
an ample divisor under the birational map $(\pi_1,\pi_2):\M_{0,5}\rightarrow \M_{0,4}\times
\M_{0,4}$, hence it is big.

\indent \textbf{3.} There is an $i\in \{1,2\}^c$ such that $D\cdot C^{1i}=0$, say $D\cdot
C^{13}=0$. Then $c_2=c_3=c_4=c_5=c>0$ and $c_1=-2c.$ Proposition \ref{m05} gives that
$D\equiv c\sum_{a,b\neq 1}\delta_{ab}=c\pi_1^*(\delta_{23}+\delta_{24}+\delta_{25})$, which
proves our contention.
\hfill $\Box$

\section{\textbf{The nef cone of $\M_{0,6}$}}
In this section we prove Theorems \ref{fulton} 
and \ref{fibrations}. The  main idea is to canonically write every divisor $D$ on $\M_{0,6}$ as a linear combination of boundary divisors with coefficients being intersection numbers with specific combinations of $F$-curves.

We first introduce a number of $1$-cycles on $\M_{0,6}$. Whenever we refer to a $1$-cycle as being a {\sl weighted sum} of $F$-curves we mean that we divide by the number of irreducible components making up the cycle. Let us fix distinct $a,b \in \{1,\ldots,6\}$. By $C_{1}^{ab}$
(respectively $C_{2}^{ab}$) we denote
 the weighted sum of  $F$-curves of type $(3:1:1:1)$ (resp.
$(2:2:1:1)$) having both points indexed by $a$ and $b$  on the spine. By $C_3^{ab}$ we denote the weighted sum of $F$-curves of type $(3:1:1:1)$ having neither $a$ nor $b$ on the spine, while $C_4^{ab}$ is the weighted sum of $F$-curves of type $(2:2:1:1)$ with $a,b$ on the same tail.

 For $a,b \in
\{2,\ldots,6\}$, let $C_{1}^{1ab}$ be the weighted sum of $F$-curves of type $(2:2:1:1)$  having
exactly one of the points indexed by elements of $\{1,a,b\}$ on the spine while the remaining two points are on one of the attached tails.  By $C_{2}^{1ab}$ (resp. $C_{3}^{1ab}$) we denote the weighted sum of $F$-curves  of
type $(3:1:1:1)$ having only one of the points (resp. two of the points) indexed by
elements of $\{1,a,b\}$ on the spine.

It may be of interest to note that  $C_{1}^{ab}$ and $C_{2}^{ab}$ are the only $G_{ab}$-invariant
$F$-curves (up to rescaling) that properly intersect $\Delta_{ab}$. Similarly
$C_{1}^{1ab}$ is the unique  $G_{1ab}$-invariant $F$-curve that properly intersects
$\Delta_{1ab}$ and  $C_{2}^{1ab}$ and $C_{3}^{1ab}$ are the only $G_{1ab}$-invariant curves of type $(3:1:1:1)$ that do not intersect $\Delta_{1ab}$ at all. 

Throughout this section we use the notation from Definition 1. To simplify things we set $\Sigma^{abc}:=\Sigma^{abc,2}_1$, $\Sigma:=\Sigma^{abc}+\delta_{abc}$ and $I+O:=I^{abc}+O^{abc}=I^{ab}+O^{ab}$, for any $a,b,c\in \{1,\ldots,6\}$.
The following lemma describes the intersections of the previously introduced curves with any divisor.

\begin{lemma}\label{intersection}

If $D\equiv  \sum_{1 \leq i \leq 6}c_i \psi_i - \sum_{ij \in \{2,\ldots,6\}}b_{1ij}\delta_{1ij}$
is any divisor on $\M_{0,6}$, then for distinct $a,b,c \in \{1,\ldots,6\}$ we have that 
\begin{multline*}
C^{ab}_{1} \cdot D = I^{ab} + \frac{1}{4}O^{ab} + \frac{1}{4}\Sigma_1^{ab,2}, \ \  
C^{ab}_{2} \cdot D = I^{ab}  - \frac{1}{3} \Sigma_2^{ab,1},\ \ 
C^{ab}_3 \cdot D=\frac{3}{4} O^{ab}+\frac{1}{4} \Sigma_1^{ab,2},\\ 
C^{ab}_4 \cdot D=\frac{1}{2} O^{ab}-\frac{1}{2} \Sigma_1^{ab,2}, \ \  \mbox{ }C^{abc}_{1} \cdot
D = \frac{1}{3}(I+O) - b_{abc} - \frac{1}{9}\Sigma^{abc}, \\ \ \
C^{abc}_{2} \cdot D = \frac{1}{3}I^{abc}+ \frac{2}{3}O^{abc} 
+ \frac{1}{9}\Sigma^{abc},\ \  \   
C^{abc}_{3} \cdot D = \frac{2}{3}
I^{abc}+ \frac{1}{3}O^{abc} + \frac{1}{9} \Sigma^{abc}.\\
\end{multline*}
\end{lemma}

\begin{proof}
This follows from standard intersection calculations as explained in for example
\cite{HMo} or \cite{Fa1}.
\end{proof}

The following sufficient criteria for a divisor on $\M_{0,6}$ to be big  will prove useful a number of times:
\begin{lemma}\label{bigdiv}
Let $\{i,j,k,l,m,n\}$ be a permutation of $\{1,\ldots,6\}$. Then any effective sum of boundary classes supported on  $\delta_{mn}$, $\delta_{il}$, $\delta_{jl}$, $\delta_{kl}$,
$\delta_{mni}$, $\delta_{mnj}$ and $\delta_{mnk}$ is big. Moreover any effective class supported on $\delta_{il},\delta_{jm},\delta_{kn}$ and on all boundaries $\delta_{abc}$ except $\delta_{ijk}$ is big as well.
\end{lemma}
\begin{proof}  
For the first statement it is enough to consider the pullback of the ample class $(\delta_{mn},\delta_{mn},\delta_{mn})$ under the birational map
$(\pi_{ij},\pi_{jk},\pi_{ik}): \M_{0,6} \rightarrow \M_{0,4} \times \M_{0,4}
\times \M_{0,4}$ whose components forget the marked points $(i,j)$, $(j,k)$ and $(i,k)$ respectively. 
To prove  the second statement  we pull back the class $(\delta_{jm}+\delta_{kn},\delta_{il}+\delta_{kn},\delta_{il}+\delta_{jm})$
via the birational map $(\pi_{il},\pi_{jm},\pi_{kn}):\M_{0,6}\rightarrow \M_{0,4}\times \M_{0,4}\times \M_{0,4}$.
\end{proof}

We have the following canonical presentation of any divisor class on $\M_{0,6}$:

\begin{prop}\label{presentation}
 Any divisor $D$ on $\M_{0,6}$ can be written as
\begin{eqnarray*} 
D
& \equiv & \sum_{a,b \in \{1, \ldots, 6\}}^{} \bigl( \frac{2}{5}  \bigl(C_{1}^{ab} \cdot D 
\bigr)  + \frac{1}{5} \bigl(C_{2}^{ab} \cdot D  \bigr) \ \bigr) \delta_{ab} +\ \\ & +  &
\sum_{a,b \in \{2, \ldots, 6\}}^{}  \Bigl (\frac{7}{10} \bigl(C_{1}^{1ab} \cdot D  \bigr)  +
\frac{1}{15} \bigl( (C_{2}^{1ab} + C_{3}^{1ab}) \cdot D  \bigr)  + \frac{4}{135}  \Sigma^{1ab}
\Bigr)  \delta_{1ab}. \\ \end{eqnarray*} 
\end{prop}

\begin{proof}
We perform two canonical modifications of the expression of any divisor on $\M_{0,6}$ in the basis referred to in Lemma \ref{basis}. In this way we get two presentations for any divisor on $\M_{0,6}$. The expression from Proposition \ref{presentation} is obtained by taking a suitable linear combination of them. Note that if $D$ is an $F$-divisor the $\delta_{ab}$ part
of the expression of $D$ is always effective. 
 
We start with a divisor $D \equiv 
\sum_{1 \leq i \leq 6}^{}c_i \psi_i - \sum_{j,k \in \{2,\ldots,6\}}^{}
b_{1jk}\delta_{1jk}$. We replace each $\psi_i$ by its average expansion provided by Lemma \ref{psi} to get that $D$ is linearly equivalent to
$$D^I=
\sum_{j,k \in \{2,\ldots,6\}}^{}\bigl( \frac{3}{10}(I+O) - b_{1jk}\bigr) \ 
\delta_{1jk}
+ \sum_{a,b \in \{1,\ldots,6\}}^{}\bigl( \frac{3}{5} I^{ab} +
 \frac{1}{10} O^{ab}) \bigr) \delta_{ab}.
$$
Next, in $D^I$ we replace each class $\delta_{1jk}$ by its average formula from Proposition \ref{average},
\begin{equation}\label{avrg}
\delta_{1jk}= \frac{2}{9} \sum_{\stackrel{a \in \{1,j,k\}}{b \in \{ 1,j,k \}^c
}}^{}\delta_{ab} - \frac{1}{3} \sum_{\stackrel{ab \in \{1,j,k\}\mbox{ }or}{ab \in
\{1,j,k\}^c}}^{}\delta_{ab} +  \frac{1}{9} \sum_{a,b \ne j,k}^{}\delta_{1ab},
\end{equation}
to get that $D$ can also be written as
$$D^{II}=
\sum_{j,k \in \{2,\ldots,6\}}^{}\bigl( \frac{3}{10}(I+O) - 
\frac{1}{9}\Sigma^{1jk} \bigr) \delta_{1jk} 
 + \sum_{a,b \in \{1,\ldots,6\}}^{}\bigl( \frac{3}{5} I^{ab} +
 \frac{1}{10}O^{ab} + \frac{1}{3} \Sigma_{1}^{ab,2}
 - \frac{2}{9} \Sigma_{2}^{ab,1}\bigr) \delta_{ab}.$$

We now write that $D\equiv \frac{3}{10}(\frac{7}{3}D^I+D^{II})$ and by using the intersection numbers computed in Lemma \ref{intersection} we get exactly the desired expression for $D$.
\end{proof}

To simplify notation we shall rewrite the expression from Proposition \ref{presentation} as

\begin{equation}\label{10/3}
\frac{10}{3}D = \frac{7}{3} D^I + D^{II} =
         \sum_{a,b \in \{1, \ldots, 6\}}^{}\zeta_{ab}\delta_{ab} +
                \sum_{a,b \in \{2, \ldots, 6\}}^{}\zeta_{1ab}\delta_{1ab}.
\end{equation}
 Thus $\zeta_{ab}=2I^{ab}+\frac{1}{3}O^{ab}+\frac{1}{3}\Sigma^{ab,2}_1-\frac{2}{9}\Sigma^{ab,1}_2$ and
$\zeta_{1ab}=I+O-\frac{20}{9}b_{1ab}-\frac{1}{9}\Sigma$. We have already seen that for an $F$-divisor all the coefficients $\zeta_{ab}$ are $\geq 0$. Moreover, in Proposition \ref{big} we prove that at most one of the coefficients $\zeta_{1ij}$ can be $<0$. If this happens, we replace $\delta_{1ij}$ by its average expression (\ref{avrg}) spreading the negativity of $\zeta_{1ij}$ among all boundary classes. We show that the resulting expression becomes effective thus proving Theorem \ref{fulton}. This procedure gives a decomposition of the nef cone of $\M_{0,6}$ into $11$ natural subcones: one described by inequalities $\zeta_{1ab}\geq 0$ for all $a,b\in \{2,\ldots,6\}$
and the remaining $10$ given by inequalities $\zeta_{1ij}\leq 0$ for $i,j\in \{2,\ldots,6\}$. More precisely we have the following result:

\begin{prop}\label{big}
Let $D$ be an $F$-divisor on
$\M_{0,6}$ with $\zeta_{1ij}<0$ for some $i,j\in \{2,\ldots,6\}$. Then $D$ is big and there exists a big effective combination of boundary classes $B_{1ij}$ such that
\begin{multline*} 
 D  \ \ \ \equiv  \ \ \   B_{1ij} \ \ + \ \ \sum_{\stackrel{a \in \{1,i,j\}}{b \in\{1,i,j\}^c }}^{} \bigl(\frac{1}{6}(C_{1}^{1ij} \cdot D) 
 + \frac{2}{3}  (C_{1}^{ab} \cdot D)\bigr) \ \delta_{ab}\ \ + \\
 + \sum_{\stackrel{a,b \in \{1,i,j\}\ \mbox{\tiny{ or }}}{a,b \in\{1,i,j\}^c }}^{} 
 \bigl(\frac{2}{5} (C_{1}^{ab}\cdot D  )  +  \frac{1}{5} 
 (C_{2}^{ab} \cdot D )\bigr) \ \delta_{ab} \ \ + \ \
 \sum_{a,b \neq i,j}^{}  \frac{2}{3} (C_{ab} \cdot D)  \ \delta_{1ab}\ ,   \\ 
\end{multline*}
 where
$C_{ab}$ is an effective sum of $F$-curves. In particular Theorem \ref{fulton} follows.
\end{prop}

\begin{proof} 
After replacing $\delta_{1ij}$ by its average expression in (\ref{10/3})
we obtain the identity (*)
$$
\frac{10}{3}D\equiv \sum_{\stackrel{a \in \{1,i,j\}}{b \in \{ 1,i,j \}^c
}}^{}(\zeta_{ab}+\frac{2}{9} \zeta_{1ij})\delta_{ab} + \sum_{\stackrel{a,b \in \{1,i,j\}\ \mbox{\tiny{ or}}}{a,b \in
\{1,i,j\}^c}}(\zeta_{ab}-\frac{1}{3}\zeta_{1ij})\delta_{ab} +  \sum_{a,b \neq i,j}(\zeta_{1ab}+\frac{1}{9}\zeta_{1ij})\delta_{1ab}.
$$

We set $\rho:=-(I+O+\Sigma)$. Since $0>\zeta_{1ij}\geq \zeta_{1ij}-\frac{5}{2} C^{1ij}_1\cdot D=-\frac{1}{6}\rho$, we obtain that $\rho>0$.
It is rather straightforward to check using Lemma \ref{intersection} that

\begin{equation}\label{2/9}
\zeta_{ab}+\frac{2}{9}\zeta_{1ij}=\frac{5}{9}C^{1ij}_1\cdot D+\frac{20}{9} C_1^{ab}\cdot D
+\frac{5}{27}\rho,\mbox{ }\mbox{ for }a\in \{1,i,j\}\mbox{ and  }b\in \{1,i,j\}^c,
\end{equation}

\begin{equation}\label{1/9}
\frac{1}{9}\zeta_{1ij}+\zeta_{1ab}=\frac{2\rho-8\zeta_{1ij}}{9}+\frac{20}{9}(I+O-b_{1ij}-b_{1ab}), \mbox{ }\mbox{ for } \{a,b\}\neq \{i,j\},
\end{equation}
while obviously $\zeta_{ab}-\frac{1}{3} \zeta_{1ij}>0$ for $a,b\in \{1,i,j\}$ or
$a,b\in \{1,i,j\}^c$.
We claim that (*) is already an effective representation of $\frac{10}{3}D$. As it turns out we can prove a little more than that.

For $a,b\in \{2,\ldots,6\}$ such that $\{a,b\}\neq \{i,j\}$ we define an effective $1$-cycle $C_{ab}$ such that $C_{ab} \cdot D=(I+O)-b_{1ij}-b_{1ab}.$ 
By passing to the complement if necessary, we may assume that $\{i,j\}\cap \{a,b\}=\emptyset$. We denote by $k$ the remaining marked point, hence $\{1,\ldots,6\}=\{1,i,j,k,a,b\}$. We then take
$C_{ab}:=\frac{1}{2}(2C^{1k}_3+C_4^{1k}+C')$, where $C'$ is the $F$-curve of type $(2:2:1:1)$ with $i,j$ and $a,b$ respectively sitting on different tails.

We define the divisor class $B_{1ij}$ by the formula
\begin{equation}
B_{1ij} = \frac{3}{10} \bigl(\frac{5\rho}{27} 
\sum_{\stackrel{a \in \{1,i,j\}}{b \in \{1,i,j\}^c}}{}\delta_{ab} -
\frac{\zeta_{1ij}}{3} \sum_{\stackrel{a,b \in \{1,i,j\}\ \mbox{\tiny{ or}}}{a,b \in \{1,i,j\}^c}}{}\delta_{ab} +
\frac{2\rho-8\zeta_{1ij}}{9} \sum_{a,b \neq i,j}{}\delta_{1ab} \bigr).
\end{equation}

All the coefficients in this expression are positive while the support of $B_{1ij}$ is $\sum_{S\neq \{1,i,j\}} \Delta_S$, which is a big divisor (use Lemma \ref{bigdiv}).
\end{proof}

We shall use Proposition \ref{big} to classify all fibrations of $\M_{0,6}$. We have already seen that an $F$-divisor $D$ for which there exists a coefficient $\zeta_{1ij}<0$, has to be big, hence it does not give rise to a fibration. The divisor $D$ is also big when $\zeta_{1ab}>0$ and $\zeta_{ab}>0$ for all $a$ and $b$ (use the existence of the ample class $\kappa_1$), so we are left with classifying nontrivial $F$-divisors $D$ for which all of the coefficients in (\ref{10/3})
are nonnegative and at least one of them is $0$. We have three cases to consider:
\begin{enumerate}
\item There are at least two coefficients $\zeta_{1ab}$ which are equal to $0$, that is, $D$ lies in the intersection of two of the subcones making up the nef cone of $\M_{0,6}$. Then we show that $D$ is the pullback of an effective divisor via a fibration $\M_{0,6}\rightarrow \M_{0,4}\times \M_{0,4}$ obtained by forgetting two disjoint pairs of marked points.
\item There is an $i\in \{1,\ldots,6\}$ such that $\zeta_{ij}=0$ for all $j\neq i$.
Then $D$ is the pullback of an effective divisor via the fibration $\pi_i:\M_{0,6}\rightarrow \M_{0,5}$ forgetting the $i$-th marked point.
\item If neither of the previous situations occurs then we show that $D$ is big.
\end{enumerate}

The following observation will come into nearly every argument in the rest of this section:

\begin{lemma}\label{the most useful one} Let $D\equiv \sum_{i=1}^6 c_i\psi_i-\sum_{i,j\in \{2,\ldots,6\}}
b_{1ij}\delta_{1ij}$ be a nontrivial $F$-divisor on $\M_{0,6}$ and $\{a,b,i,j,m,n\}$ a permutation of
$\{1,\ldots,6\}$. If $\zeta_{ij}=\zeta_{ab}=0$  then $\zeta_{mn}>0$.
Moreover, if $\zeta_{ia}=\zeta_{ib}=\zeta_{im}=\zeta_{in}=0$ then $\zeta_{ij}=0$ as well. \end{lemma}

\begin{proof} Without loss of generality we may assume that $\zeta_{12}=\zeta_{34}=\zeta_{56}=0$ and we prove that in this case $D$ is trivial. Our assumption implies that $C^{12}_1\cdot D=0$,
from which we can write that $c_i+b_{12i}=-(c_1+c_2)$, for all $i\in \{1,2\}^c$.
Similarly $c_j+b_{34j}=-(c_3+c_4)$, for all $j\in \{3,4\}^c$ and $c_k+b_{56k}=-(c_5+c_6)$ when $k\in \{5,6\}^c$. It is easy to see that these
relations imply that all the $c_i$'s are equal, that is, $c_i=c$ for $i\in \{1,\ldots,6\}$, hence $b_{12i}=-3c$, for each $i\in \{1,2\}^c$. Similarly $b_{34j}=-3c$ for $j\in \{3,4\}^c$ and $b_{56k}=-3c$ for each $k\in \{5,6\}^c$. 

On the other hand $C^{12}_2\cdot D=0$ which implies in particular that $c_1+c_2-b_{134}-b_{156}=0$, thus  giving that $c_i=0$ for all $i\in \{1,\ldots,6\}$. It immediately follows that the boundary coefficients must vanish too, hence $D$ is trivial.

For the second part, let us assume that $\zeta_{12}=\cdots =\zeta_{15}=0$ and we prove that $\zeta_{16}=0$. Since $C^{1i}_1\cdot D=0$ for all $i\in \{2,\ldots,5\}$ we have that $b_{1i6}=-I^{1i6}$ for all $i\in \{2,\ldots,5\}$ which yields $C^{16}_1\cdot D=0$. We also know that $C^{1i}_2\cdot D=0$ which turns out to be equivalent with $2c_1+I+O=0$. It follows that $C^{16}_2\cdot D=0$ as well, hence $\zeta_{16}=0$.
\end{proof}
We proceed with the classification of all $F$-divisors on $\M_{0,6}$. The next lemmas deal with the first two situations:
\begin{lemma}\label{type1}
Let $D\equiv \sum_{i=1}^6 c_i\psi_i-\sum_{i,j\in \{2,\ldots,6\}} b_{1ij}\delta_{1ij}$ be an $F$-divisor on $\M_{0,6}$ such that two coefficients $\zeta_{1ij}$ vanish, say $\zeta_{1ij}=\zeta_{1kl}=0$, where $\{1,i,j,k,l,m\}$ is a permutation of $\{1,\ldots,6\}$. Then $D$ is the pullback of an effective divisor via the fibration $\phi=(\pi_{ij},\pi_{kl}):\M_{0,6}\rightarrow \M_{0,4}\times \M_{0,4}$ whose components forget the marked points labelled $(i,j)$ and
$(k,l)$ respectively.
\end{lemma}
\begin{proof} We use the notation from the proof of Proposition \ref{big}. Without loss of generality we may assume that $\zeta_{123}=\zeta_{145}=0.$ From
(\ref{1/9}) we obtain that $\rho=-(I+O+\Sigma)=0$ and that $C_{16}\cdot D=0$ which implies that $C^{16}_3\cdot D=C^{16}_4\cdot D=0$. Since $\zeta_{123}=0$, we can also write that $\rho/15=C^{123}_1\cdot D=0$ and similarly $C^{145}_1\cdot D=0$.
Thus the intersection numbers of $D$ with every component of $C^{16}_3,C^{16}_4,C^{123}_1$ and $C^{145}_1$ respectively has to be $0$. This gives rise to $28$ relations between the coefficients $c_i$ and $b_{1ij}$. By writing out these relations it turns out that (\ref{10/3}) can be rewritten as
\vskip 4pt
\noindent
$D\equiv (\alpha+\beta)(\delta_{124}+\delta_{134}+\delta_{125}+\delta_{135})+\alpha (\delta_{12}+\delta_{13}+\delta_{62}+\delta_{63}+\delta_{146}+\delta_{156}+\delta_{16}+\delta_{23})+$
\newline
\noindent
$+\beta (\delta_{14}+\delta_{15}+\delta_{64}+\delta_{65}+\delta_{126}+\delta_{136}+\delta_{16}+\delta_{45}), $
\vskip 4pt
\noindent where $\alpha=\frac{2}{3}(c_1+c_2)\geq 0$ and $\beta=\frac{2}{3}(c_1-c_2)\geq 0$. In order to finish the proof it is enough to notice that $D\equiv \alpha\ \pi_{45}^*(\delta_{12}+\delta_{23}+\delta_{13})+\beta \ \pi_{23}^*(\delta_{14}+\delta_{45}+\delta_{15})$.
\end{proof}

\noindent \textbf{Remark.} Lemma \ref{type1} also shows that if $D$ is a nontrivial $F$-divisor on $\M_{0,6}$ then at most two of the coefficients $\zeta_{1ab}$ can vanish.

\begin{lemma}\label{type2}
Let $D$ be an $F$-divisor such that there exists $i\in \{1,\ldots,6\}$ with $\zeta_{ij}=0$ for all $j\neq i$. Then $D$ is the pullback of an effective divisor under the projection $\pi_i:\M_{0,6}\rightarrow \M_{0,5}$ dropping the $i$-th marked point.
\end{lemma}
\begin{proof} Clearly we can assume that $i=1$. The hypothesis $\zeta_{1j}=0$ for $j\in \{2,\ldots,6\}$ is equivalent to $C^{1j}_1\cdot D=C^{1j}_2\cdot D=0$ for all $j\in \{2,\ldots,6\}$. This gives that $b_{1ij}=-I^{1ij}$ for all $i,j\in \{2,\ldots,6\}$. It also follows that $2c_1+I+O=0$ and $I+O+\Sigma=0$.
Then (\ref{10/3}) reads $$D\equiv \sum_{i,j\in \{2,\ldots,6\}} \frac{2}{3}I^{ij}(\delta_{ij}+\delta_{1ij})=\pi_1^*(\sum_{i,j\in \{2,\ldots,6\}}
\frac{2}{3}I^{ij}\delta_{ij}).$$
\end{proof}
 
\vskip 3pt
\noindent
\emph{Proof of Theorem \ref{fibrations}, Part 2}. We start with a nontrivial $F$-divisor $D$ on $\M_{0,6}$ with $\zeta_{1ab}\geq 0$ for all $a,b\in \{2,\ldots,6\}$ and such that at most one coefficient $\zeta_{1ab}$ is equal to $0$, say $\zeta_{123}=0$. Moreover, we can assume that for each $i$ there is a $j\neq i$ such that $\zeta_{ij}\neq 0$. Then we show that $D$ is big. Note that the case $\zeta_{1ab}>0$ for all $a,b\in \{2,\ldots,6\}$ is similar (and simpler).

Lemma \ref{the most useful one} limits the number of coefficients $\zeta_{ij}$ that can vanish and a case by case analysis shows that we can always find sufficiently many boundaries $\delta_{ij}$ on which $D$ is supported. Then we apply Lemma \ref{bigdiv} to conclude that $D$ is big.
\hfill $\Box$
\vskip 3pt

\section{\textbf{The Mori cone of $\M_{g}$}}

In this section we show that $\NE_{1}(\M_g)$ is spanned by $F$-curves for all $g\leq 13$. To prove this, it is enough to show that every $S_g$-invariant extremal ray on $\M_{0,g}$ is generated by an $F$-curve (cf. \cite{GKM}, Theorem 0.3). We achieve this inductively by writing every nontrivial $S_g$-invariant nef divisor on $\M_{0,g}$ as a sum $K_{\M_{0,g}}+\sum_{S} a_S\delta_S$, where $0\leq a_S\leq 1$ for all $S$. We also notice that for any $g\geq 14$ there are $S_g$-invariant $F$-divisors on $\M_{0,g}$ not of this form, thus hinting that the nature of $\M_{g}$ changes in a subtle way when $g=14$. Finally we present a combinatorial set-up enabling us to compute Mori cones of moduli spaces of $1$ and $2$-pointed curves of genus $g\leq 8$. It is clear that in the same way at least a couple of other cases can be settled as well.

We start by setting some notation. We denote by $\widetilde{M}_{0,n}:=\M_{0,n}/S_n$ and we identify divisors on $\widetilde{M}_{0,n}$ with $S_n$-invariant divisors on $\M_{0,n}$. The spaces $\widetilde{M}_{0,n}$ are interesting for their own sake. For instance $\widetilde{M}_{0,2g+2}$ is isomorphic to the closure in $\M_g$ of the locus of hyperelliptic curves of genus $g$. For $2\leq i<\lfloor n/2 \rfloor$, we set
$B_i := \sum_{ \stackrel{S \subset \{1,\ldots,g\}}{|\S| = i}}^{}\delta_{\S}.$
When $i=n/2$ we define $B_i:=\sum_{\stackrel{S \subset \{1,\ldots,g\}}{|S|=i,1\in S}} \delta_S.$ 

 Keel and McKernan proved the following results about the Mori theory of $\widetilde{M}_{0,n}$ (cf. \cite{KMcK}, Theorem 1.3):
\begin{prop}\label{mckernan}
\begin{enumerate}
\item
The effective cone $\overline{NE}^1(\widetilde{M}_{0,n})$ is generated by the classes of the divisors $B_i$ for $2\leq i\leq \lfloor n/2 \rfloor$. Any nontrivial nef divisor on $\widetilde{M}_{0,n}$ is big.
\item For $n\leq 7$ the cone of curves $\overline{NE}_1(\M_{0,n})$ is generated by $F$-curves. We also have that for  $n\leq 11$ the cone $\overline{NE}_1(\widetilde{M}_{0,n})$ 
is spanned by $F$-curves.
\end{enumerate}
\end{prop}

\noindent \textbf{Remark.} The previous result combined with Theorem 0.3 from \cite{GKM} gives that $\NE_1(\M_{g,n})$ is spanned by $F$-curves whenever $g+n\leq 7$. We also obtain that for $g\leq 11$ every $F$-divisor on $\M_g$ is
nef. We shall extend this result for all $g\leq 13$. 
\vskip 7pt
We recall that for any $S \subset \{1,\ldots,n\}$ such
that $|S|,|S^c|\geq 2$, there is an isomorphism 
$$
\phi:\M_{0,|\S|+1} \times \M_{0,|\S^c|+1} \longrightarrow \Delta_{\S} \subseteq
\M_{0,n}
$$
given by attaching a rational $(|S|+1)$-pointed curve to 
a rational $(|S^c|+1)$-pointed curve at a point $x$.
It turns out that $\phi$ induces an isomorphisms between Mori cones
$\NE_1(\Delta_S)=\NE_1(\M_{0,|S|+1})\times \NE_1(\M_{0,|S^c|+1})$
(cf. \cite{KMcK}, Lemma 3.8). Moreover, if $\pi_1:\Delta_S\rightarrow \M_{0,|S|+1}$ and $\pi_2:\Delta_S\rightarrow \M_{0,|S^c|+1}$ are the two projections then 
$$N_{\Delta_S/\M_{0,n}}\equiv (\pi_1)^*(-\psi_x)+ (\pi_2)^*(-\psi_x)$$
(cf. \cite{KMcK}, Lemma 4.5).
Since the tautological classes $\psi_x$ are nef (cf. \cite{Kap}), it follows that  $\Delta_{S}$ has {\sl anti-nef} normal bundle, that is, $C\cdot \Delta_S\leq 0$ for every irreducible curve $C\subseteq \Delta_S$. We shall often use certain maps between moduli spaces which we call {\sl boundary restrictions}:

\begin{definition} For $m,n\geq 3$ and $ n_{x_1},\ldots,n_{x_m}\geq 1$ such that $n=n_{x_1}+\cdots+n_{x_m}$, we define the map
$
\nu : \M_{0,m} \longrightarrow \M_{0,n}
$ which takes a rational $m$-pointed curve $(C,x_1,\ldots,x_m)$ to 
a rational $n$-pointed curve by attaching a fixed rational $(n_{x_i}+1)$-pointed curve at each point $x_i$ for $1\leq i\leq m$.  
\end{definition}

Note that if
$n_{x_i} =1$ for some $i$ then this amounts to not having attached anything at
$x_i$. Moreover any composition of 
boundary restrictions will be homotopic to a single boundary restriction, in particular
they will induce the same map in homology. We make the following simple observation:

\begin{prop}\label{iteration}
Given integers $g,n$ with $g+n\geq 8$, to conclude that $\NE_1(\M_{g,n})$ is generated by $F$-curves it suffices to prove that for all $F$-divisors $D$ on
$\M_{0,g+n}/S_g$ and for all boundary restrictions $\nu:\M_{0,k}\rightarrow \M_{0,g+n}$, where $8\leq k\leq g+n$, the pullback  $\nu^*(D)$ is a nonnegative combination of boundary divisors.
\end{prop}
\begin{proof}
We apply \cite{GKM}, Theorem 0.3. We start with an $F$-divisor $D$ on $\M_{0,g+n}/S_g$ and we want to 
show that $D$ is nef. Since $D$ is a nonnegative combination of boundary divisors  we only have to show that $C\cdot D\geq 0$ for all irreducible curves in a boundary divisor $\Delta_S\cong \M_{0,|S|+1}\times \M_{0,|S^c|+1}$. By hypothesis $D_{|\M_{0,|S|+1}}$ 
and $D_{|\M_{0,|S^c|+1}}$ are both effective combinations of boundary, hence we have to test the
 nefness of $D$ only against curves sitting in the boundary of $\M_{0,|S|+1}$ and of $\M_{0,|S^c|+1}$ 
 and we can descend all the way to a moduli space $\M_{0,n}$ with $n\leq 7$. Since in this range the 
 $F$-curves generate $\overline{NE}_1(\M_{0,n})$ (cf. Proposition \ref{mckernan}), the conclusion follows. 
\end{proof}

Let us consider a boundary restriction $\nu:\M_{0,m}\rightarrow \M_{0,n}$ given by the partition 
$(n_{x_1},\ldots,n_{x_m})$ of $n$ where we assume that $n_{x_j}\geq 2 \Leftrightarrow j\in \{1,\ldots,r\}$. 
We also denote by $A:=\{x_{r+1},\ldots,x_m\}$ the set of remaining marked points, hence $n_y=1$ for all $y\in A$.
For $2\leq i\leq m-2$ and for $S\subset \{x_1,\ldots,x_r\}$ we define

$$B_i^S:=\sum_{\stackrel{T\subset \{1,\ldots,m\},|T|=i}{T\cap\{x_1,\ldots,x_r\}=S}} \delta_T.$$

The adjunction formula for a bounday restriction $\nu:\M_{0,m}\rightarrow \M_{0,n}$ reads 
\begin{equation}\label{adjunction} 
\nu^*(K_{\M_{0,n}})=K_{\M_{0,m}}+\sum_{n_x\geq 2}\psi_x.
\end{equation}
The next statement describes the effect a boundary restriction has on homology.

\begin{prop}\label{restriction}
 Let $D \equiv \sum_{i=2}^{\lfloor n/2 \rfloor}r_i B_i$ be a divisor on $\widetilde{M}_{0,n}$. If
 $\nu:\M_{0,m}\rightarrow \M_{0,n}$ is a boundary restriction then
$
\nu^*D \equiv  \sum_{\stackrel{S \subset \{x_1,\ldots,x_r\}}{|S|\leq \lfloor r/2 \rfloor, i\geq |S|}}^{} c_i^S B_i^S,
$
where 
$$
c_i^S:=r_{i+\sum_{x \in S} n_x -|S|} - 
\frac{(m-i)(m-1-i)(\sum_{x\in S} r_{n_x})+i(i-1)(\sum_{x\in S^c}r_{n_x})}{(m-1)(m-2)}.
$$
\end{prop}

\begin{proof} We start with the case when only one $n_x$ is $\geq 2$. We obtain that $\nu^*(B_i)=B_{i}^{\emptyset}+B_{i-n_x+1}^{x}-(?i)\psi_i$, where $(?i)=1$ if $i=n_x$ and $0$ otherwise. By iteration, in the case when $n_{x_j}\geq 2$ for $1\leq j\leq r$, we can write that
$$\nu^*(B_i)=\sum_{S\subset \{x_1,\ldots,x_r\}} B_{i-\sum_{x\in S}n_x+|S|}^S-\sum_{\stackrel{x\in \{x_1,\ldots,x_r\}}{n_x=i}}\psi_{x}.$$
To read this formula correctly, when $i=n/2$ the first sum is taken only over the subsets $S\subset \{x_1,\ldots,x_r\}$ containing $x_1$, that is, we do not count both $S$ and $S^c$. Moreover we make the convention that $B_j^S=0$ whenever $j<|S|$ or $j\geq m-1$.
Now replacing each $\psi_x$ by its average formula from Lemma \ref{psi} we get the desired formula.
\end{proof}

We have already seen that $\NE_1(\M_{g,n})$ is spanned by $F$-curves whenever $g+n\leq 7$. We now present an entirely combinatorial computation of $\NE_1(\M_{6,2})$ which can serve as a model for other cases when $g+n$ is rather small.

\begin{prop}\label{m62}
Every $F$-divisor on $\M_{0,8}/S_6$ is linearly equivalent to an effective combination of boundary divisors. It follows that $\NE_1(\M_{6,2})$ is generated
by $F$-curves.
\end{prop}

\begin{proof}
Let us denote by $X:=\M_{0,8}/S_6$ and by $x,y\in \{1,\ldots,8\}$ the marked points on which $S_6$ fails to act. We identify divisors on $X$ with $S_6$-invariant divisors on $\M_{0,8}$. A basis for $\mbox{Pic}(X)$ is given by the following collection of divisor classes:

\vskip 5pt
\noindent
$\delta_{x1}:=\sum_{a\neq x,y} \delta_{xa},\mbox{ } \delta_{y1}:=\sum_{a\neq x,y} \delta_{ya},\mbox{ } \delta_{x2}:=\sum_{a,b\in \{x,y\}^c} \delta_{xab}, \mbox{ }\delta_{y2}:=\sum_{a,b\in \{x,y\}^c} \delta_{yab},$
\vskip 7pt
\noindent 
$\delta_2:=\sum_{a,b\in \{x,y\}^c} \delta_{ab}, \mbox{ }\delta_3:=\sum_{a,b,c\in \{x,y\}^c} \delta_{abc}, \mbox{ } \delta_{xy1}:=\sum _{a\neq x,y} \delta_{axy},$
\vskip 7pt
\noindent
$\delta_{xy2}=\delta_4:=\sum_{a,b\in \{x,y\}^c} \delta_{xyab}\mbox { and }\delta_{x3}=\delta_{y3}:=\sum_{a,b,c\in \{x,y\}^c} \delta_{xabc}.$

There is a $10$-th $S_6$-invariant divisor class on $\M_{0,6}$, namely $\delta_{xy}$, which can be expressed in this basis using our average formula from Proposition \ref{average}:
\begin{equation}\label{avm8}
\delta_{xy}=\frac{1}{6}(\delta_{x1}+\delta_{y1})+\frac{4}{15}(\delta_{x2}+\delta_{y2})+\frac{3}{20} \delta_{x3}-\frac{1}{5}\delta_3-\frac{2}{3} \delta_{xy1}-\frac{2}{5} \delta_{xy2}-\frac{1}{15} \delta_2.
\end{equation}
 We now start with an arbitrary $F$-divisor on $X$: 
\begin{equation}\label{d1}
D\equiv b_{x1}\delta_{x1}+b_{y1}\delta_{y1}+b_{x2}\delta_{x2}+b_{y2}\delta_{y2}+b_{x3}\delta_{x3}+b_{xy1}\delta_{xy1}+b_{xy2}\delta_{xy2}+b_2\delta_2+b_3\delta_3.
\end{equation}

The coefficients of $D$ are subject to $28$ $F$-inequalities coming from all equivalence classes of partitions of $\{1,\ldots,8\}$ in four subsets modulo the
$S_6$ action. By $(3_x,2_y,2,1)$ for instance we shall denote a partition of type $(3:2:2:1)$ such that $x$ is contained in the subset with $3$ elements and $y$ is part of one the subsets with $2$ elements. The proof that $D$ is linearly equivalent to an effective boundary consist of two parts and is similar to the proof of Theorem \ref{fulton}:

\noindent \textbf{(i)} If $b_2\geq 0$ then then all coefficients in (\ref{d1}) are nonnegative.

\noindent \textbf{(ii)} If $b_2<0$ then we substitute $\delta_2$ using (\ref{avm8})
and we call the resulting divisor class $D_2\equiv D$. We then take the linear combination $D':=5D+D_2\equiv 6D$ and we show that the resulting expression, 
\begin{equation}\label{d3}
D'\equiv (6b_{x1}+\frac{5}{2}b_2)\delta_{x1}+(6b_{y1}+\frac{5}{2} b_2)\delta_{y1}+(6b_{x2}+4b_2)\delta_{x2}+(6b_{y2}+\delta_2)\delta_{y2}+
\end{equation}

\noindent $+(6b_{x3}+\frac{9}{4}b_2)\delta_{x3}+(6b_{xy1}-10b_2)\delta_{xy1}+
(6b_{xy2}-6b_2)\delta_{xy2}+(6b_3-3b_2)\delta_3-15b_2 \delta_{xy},$
\vskip 7pt 
\noindent
is effective. We present step (ii), step (i) being similar, only simpler. 

We thus assume that $b_2<0$ and we show that all coefficients in (\ref{d3}) are nonnegative. We start with the $\delta_{xy2}$ coefficient whose nonnegativity 
follows from the $F$-inequality corresponding to the partition $(2_{xy},2,2,2).$
The fact that the $\delta_3$ coefficient is $\geq 0$ comes using $(2_{xy},2,2,2)$ and $(4_{xy},2,1,1)$. 

The other inequalities are slightly more involved. We first prove that $b_{x1}\geq 0$ by combining $(5_y,1_x,1,1), (4_y,2,1_x,1)$ and $(3,3,1_x,1_y)$.
By symmetry we also obtain that $b_{y1}\geq 0$.

Next, by adding together $(4_x,2_y,1,1),(4_y,2_x,1,1), (3,3,1_x,1_y), (3_x,2_y,2,1)$ and $(3_y,2_x,2,1)$, we find that
$$4(b_{x2}+b_{y2})+2b_2\geq 3(b_{x1}+b_{y1})\geq 0,$$
and since we have assumed $b_2<0$ we get that $b_{x2}+b_{y2}\geq 0$.
Now $(3_x,3_y,1,1)$ gives that $2b_{x3}+b_2\geq b_{x2}+b_{y2}\geq 0$, thus we also have that $b_{x3}\geq 0$. We can now prove that the remaining coefficients in (\ref{d3}) are nonnegative as well.

We start with the $\delta_{x3}$ coefficient, which is nonnegative because $2b_{x3}+b_2\geq 0$ and $b_{x3}\geq 0$. To deal with the $\delta_{x2}$ coefficient we combine $(4_y,2_x,1,1), (5_y,1_x,1,1)$ with $b_{x3}\geq 0$.
By symmetry, the $\delta_{y2}$ coefficient is also $\geq 0$. For the $\delta_{x1}$ coefficient we use that $\delta_{x1}\geq 0$ together with $(5_y,1_x,1,1)$. Again, by symmetry, the $\delta_{y1}$ coefficient is also $\geq 0$.

We are left with the $\delta_{xy1}$ coefficient whose nonnegativity follows from $b_2<0$ together with $(4,2_{xy},1,1)$ and $(2_{xy},2,2,2)$. Note that we only used $10$ of the $28$ $F$-inequalities.
\end{proof}

We use Proposition \ref{m62} to compute the Mori cones of $\M_{g,1}$ when $g\leq 8$:

\begin{prop}
The cone $\NE_1(\M_{g,1})$ is generated by $F$-curves for all $g\leq 8$.
\end{prop}

\begin{proof}
Since the case $g\leq 6$ is settled by Proposition \ref{mckernan} we only need to deal with $\M_{7,1}$ and $\M_{8,1}$. We only present the $g=8$ case, $g=7$ being similar. According to Proposition \ref{iteration}, it suffices to show that for any $F$-divisor $D$ on $\M_{0,9}/S_8$ we have that (i) D is linearly equivalent to an effective sum of boundaries, and (ii) for a boundary restriction $\nu:\M_{0,8}\rightarrow \M_{0,9}$, $\nu^*(D)$ is linearly equivalent to an effective sum of boundaries.

To prove (i) we denote by $Y:=\M_{0,n}/S_{n-1}$ and by $x$ the marked point on which $S_{n-1}$ does not act. A basis for $\mbox{Pic}(Y)$ is given by the classes $\delta^{\{x\},1}_i$ for $i=1,\ldots,n-3$. In this case the $S_{n-1}$-invariant boundary classes on $\M_{0,n}$ are independent which considerably reduces the combinatorial complexity of the problem. We write (uniquely) the class of any $F$-divisor $D$ on $Y$ in this basis and the positivity of the coefficients follows in a straightforward way from the $F$-inequalities. We
 omit the details.

For (ii) it is enough to nottice that for any boundary restriction $\nu:\M_{0,8}\rightarrow \M_{0,9}$ the pullback $\nu^*(D)$ is an $F$-divisor on $\M_{0,8}/S_6$ so by Proposition \ref{m62} it is equivalent to an effective sum of boundary classes.
\end{proof}

We are now going to prove that $\overline{NE}_1(\widetilde{M}_{0,n})$ is spanned by $F$-curves for all $n\leq 13$. We use our Theorem \ref{fulton} to give a Mori theoretic sufficient condition for an extremal ray on $\M_{0,n}$ to be generated by an $F$-curve. The next theorem is an improvement of \cite{KMcK}, Theorem 1.2.
We recall that $\Delta$ denotes the total boundary in $\M_{0,n}$.
\begin{theorem}\label{mori}
Let $R$ be an extremal ray in $\overline{NE}_1(\M_{0,n})$. If there exists a  nonempty effective $\mathbb Q$-divisor $G$ on $\M_{0,n}$ such that $\Delta-G$ is also effective and $(K_{\M_{0,n}}+G)\cdot R\leq 0$, then $R$ is contractible and it is spanned by an $F$-curve.
\end{theorem}

\begin{proof} We follow the same lines as in \cite{KMcK}. Let us write $G\equiv \sum_{S} a_S \Delta_S$, with $0\leq a_S\leq 1$. 
We claim that the ray $R$ descends to some boundary divisor. Suppose this is not the case, hence $R\cdot \Delta_S\geq 0$ for all $S$ and $R\cdot K_{\M_{0,n}}\leq 0$.
Since there exists an ample divisor on $\M_{0,n}$ having the same support as $\Delta$, namely the tautological divisor $\kappa_1$, we find that $R$ is generated by a contractible curve $C$. We denote by $f:\M_{0,n}\rightarrow Y$ the contraction. The curve $C$ does not come from the boundary hence $f_{|\Delta}$ is finite and we can apply \cite{KMcK}, Proposition 2.5 to conclude that the exceptional locus $\mbox{Exc}(f)$ is a curve. On the other hand, using the deformation theoretic bound for the dimension of the Hilbert scheme (cf. \cite{K}, Theorem 1.14)  $$\mbox{dim}_{[C]}\mbox{Hilb}(\M_{0,n})\geq -K_{\M_{0,n}}\cdot C+n-6\geq 1\mbox{ }\mbox{  (for }n\geq 7),$$
we conclude that $C$ deforms inside $\M_{0,n}$ which contradicts that $\mbox{Exc(f)}$ is a curve. This argument breaks down for $n\leq 6$ but in that case we can invoke directly Theorem \ref{fulton} and finish the proof. Thus in any case we may assume that $R$ is contained in some boundary divisor $\Delta_T$ and since 
$\overline{NE}_1(\Delta_T)=\overline{NE}_1(\M_{0,|T|+1})\times \overline{NE}_1(\M_{0,|T^c|+1})$ we may as well assume that say, $R\subseteq \overline{NE}_1(\M_{0,|T|+1}).$

If $m:=|T|+1$, we denote by $\nu:\M_{0,m}\rightarrow \M_{0,n}$ the corresponding boundary restriction and by $x$ the point of attachment of the fixed $(|T^c|+1)$-pointed rational curve. We replace $G$ by the effective divisor 
$G':=G+(1-a_T)\Delta_T$. The boundary $\Delta_T$ has anti-nef normal bundle 
 hence $(K_{\M_{0,n}}+G')\cdot R\leq 0$.

By adjunction, $\nu^*(K_{\M_{0,n}})=K_{\M_{0,m}}+\psi_x$, while according to
Proposition \ref{restriction} we have that $\nu^*(G')=\tilde{G}-\psi_x$, where $\tilde{G}$ is an effective divisor such that $\Delta_{\M_{0,m}}-\tilde{G}$ is  effective too. Thus $(K_{\M_{0,m}}+\tilde{G})\cdot R\leq 0$, that is, we have exactly the initial situation on a lower dimensional moduli space and the conclusion follows inductively. 
\end{proof}

Now we show that for $n\leq 13$ every extremal ray on $\widetilde{M}_{0,n}$ satisfies the conditions from Theorem \ref{mori}.  We start with an extremal ray $R\subseteq \overline{NE}_1(\widetilde{M}_{0,n})$ and denote by  $E$ a supporting nef divisor of $R$. Proposition \ref{mckernan} gives that $E$ is big, that is, $E\in \mbox{int}(\overline{NE}^1(\M_{0,n}))$.  Since 
$
K_{\M_{0,n}} = \sum_{j=2}^{\lfloor\frac{n}{2}\rfloor} \bigl(j(n-j)/(n-1)-2\bigr) B_j
$,
clearly $-K_{\M_{0,n}}$ is not effective for $n\geq 7$. Following Keel and McKernan we intersect the line segment in $NS(\M_{0,n})$ joining $-K_{\M_{0,n}}$ and  
 $E$ with the boundary of 
$\overline{NE}^1(\widetilde{M}_{0,n})$ to get a symmetric boundary class
$\Delta_{E}$ such that $\lambda E\equiv K_{\M_{0,n}}+\Delta_E$ for some $\lambda >0$.  We can write  
$\Delta_E \equiv \sum_{i = 2}^{\lfloor \frac{n}{2} \rfloor} r_i B_i
$, where $r_i\geq 0$ (cf. Proposition \ref{mckernan}).  That $\Delta_E$ is on an extremal face of the cone means $r_i = 0$ for at least one $i$ with $2\leq i\leq \lfloor n/2 \rfloor$. If we can prove that $r_i\leq 1$ for all $i$, then Theorem \ref{mori} gives that $R$ is generated by an $F$-curve. 

To achieve this we write out all $F$-inequalities for the nef divisor $K_{\M_{0,n}}+\Delta_E$: We define the function $f(a,b,c,d)$ to be $2$ minus the number of variables equal to $1$. For any partition $(a,b,c,d)$ of $n$ into 
positive integers we consider the associated $F$-curve given by a boundary restriction $\nu:\M_{0,4}\rightarrow \M_{0,n}$. Then using (\ref{adjunction}), 
$$
(K_{\M_{0,n}}+\Delta_E)\cdot \nu(\M_{0,4})=f(a,b,c,d)+r_{a+b}+r_{a+c}+r_{a+d}-r_a-r_b-r_c-r_d\geq 0.
$$

\begin{theorem}
For $n\leq 13$ any nontrivial nef divisor on $\widetilde{M}_{0,n}$ is of the form $K_{\M_{0,n}}+\Delta_E$, with $0\leq \Delta_E\leq \Delta$. It follows that $\NE_1(\M_g)$ is generated by $F$-curves for all $g\leq 13$.
\end{theorem}

\begin{proof} We start with the nef divisor $K_{\M_{0,n}}+\Delta_E$, where $\Delta_E\equiv \sum_{j=2}^{\lfloor n/2 \rfloor} r_j B_j$, where $r_j\geq 0$ for all $j$ and there is $2\leq i\leq \lfloor n/2 \rfloor$ such that $r_i=0$. By using all $F$-inequalities the coefficients $r_j$ are subject to, we conclude that $r_j\leq 1$. We carry this out only for $n=13$, the case $n=12$ being entirely similar. We list all $F$-inequalities for $n=13$:
\vskip 4pt
\noindent $(1)\mbox{ } 3r_2\geq r_3+1,\mbox{ }\mbox{ }(2)\mbox{ }2r_3\geq r_4, \mbox{ }\mbox{ }(3) \mbox{ }r_2+2r_4\geq r_3+r_5,\mbox{ }\mbox{ }(4) \mbox{ }r_2+2r_5\geq r_4+r_6,$
\newline
\noindent $(5)\mbox{ }r_2+r_6\geq r_5,\mbox{  }\mbox{ }(6)\mbox{ }1+2r_3+r_4\geq 2r_2+r_5, \mbox{ }\mbox{ }(7) \mbox{ }1+r_4+r_5\geq r_2+r_6,$
\newline
\noindent $(8)\mbox{ }1+r_3+r_5\geq r_2+r_4,\mbox{ }\mbox{ }(9)\mbox{ }1+r_3+2r_6\geq r_2+2r_5, \mbox{  } \mbox{ }(10) \mbox{ }1+2r_4\geq 2r_3,$
\newline
\noindent $(11)\mbox{ }1+r_6\geq r_3, \mbox{ }\mbox{ }(12)\mbox{ } 1+3r_5\geq 3r_4, \mbox{  }\mbox{ }(13)\mbox{ } 2+3r_4\geq 3r_2+r_6,$
\newline
\noindent
$(14) \mbox{ } 2+r_4+2r_5\geq 2r_2+r_3+r_6, \mbox{ }\mbox{ } (15)\mbox{ }2+2r_6\geq 2r_2+r_5,\mbox{ }\mbox{ }(16)\mbox{ } 2+r_5+r_6\geq r_2+2r_3,$
\newline
\noindent
$(17) \mbox{ } 2+r_5+2r_6\geq r_2+r_3+2r_4,\mbox{ and finally }(18)\mbox{ } 2+3r_6\geq 3r_3+r_4$.
\vskip 5pt
From (1) we see that $r_2>0$. We have four cases:

\noindent \textbf{(i) $r_3=0$.} Then from (2) we have $r_4=0$ while from (6) we get $r_5< 1$ and $r_2\leq 1/2$. Assume now that $r_6\geq 1$. Then (3) and (7) combined give $r_5=r_2$ and $r_6=1$. From (13) we get that $r_2\leq 1/3$ while (4) gives that $r_2\geq 1/3$, hence $r_2=r_5=1/3$. Thus either $r_j<1$ for all $j$ or else $\Delta_E\equiv \frac{1}{3}(B_2+B_5)+B_6$.
\vskip 4pt
\noindent \textbf{(ii) $r_4=0$.} Use (10), (13) and (3) to get that $r_3\leq 1/2,r_2\leq 2/3$ and $r_5\leq 2/3$. We assume again that $r_6\geq 1$. Then (3) and (7) are compatible only when $r_2=r_5$ and $r_3=0$, that is, we are back to case (i).
\vskip 4pt
\noindent \textbf{(iii) $r_5=0$.} Inequality (12) yields $r_4\leq 1/3$ while from (13) we get that $r_2\leq 1$. Moreover if $r_2=1$ then $r_6=0,r_4=1/3$ and $
r_3=1/3$, so in this way we get our second exceptional case, $\Delta_E\equiv B_2+\frac{1}{3}(B_3+B_4)$. On the other hand if $r_2<1$ then from (4) we have that $r_6<1$ and finally from (10) we obtain that $r_3\leq 5/6$.
\vskip 4pt
\noindent \textbf{(iv) $r_6=0$.} From (18) we have that $r_3\leq 2/3$ while (5) and (15) give that $r_5\leq 2/3$ and $r_2\leq 1$. Moreover $r_2=1$ implies $r_5=0$ so we are back to case (iii). Then we  use (12) which gives $r_4\leq 1$. If $r_4=1$ from (18) we have $r_3\leq 1/3$ while from (2) $r_3\geq 1/2$, a contradiction, so this last case does not occur.
\end{proof}

\noindent \textbf{Remark.} For $n\leq 11$ every nontrivial $F$-class on $\widetilde{M}_{0,n}$ is of the form $K_{\M_{0,n}}+\Delta_E$, where $\Delta_E$ is a {\sl pure boundary}, that is, $\Delta_E=\sum_{i} r_iB_i$ where $0\leq r_i<1$ (cf. \cite{KMcK}, Corollary 5.3). The previous proof shows that on $\widetilde{M}_{0,13}$ there are exactly two $F$-classes not of this form:
$$K_{\M_{0,13}}+\frac{1}{3}(B_2+B_5)+B_6\mbox{ and } K_{\M_{0,13}}+B_2+\frac{1}{3} (B_3+B_4),$$
On $\widetilde{M}_{0,12}$ there is just one such class, namely $K_{\M_{0,12}}+(B_2+B_5)/3+B_6$. For $n\geq 14$ it is no longer true that any nontrivial $F$-
divisor on $\widetilde{M}_{0,n}$ is numerically equivalent to $K_{\M_{0,n}}+\Delta_E$, where $\Delta_E=\sum_i r_iB_i$ with $0\leq r_i\leq 1$. When $n=14$ for instance, the $F$-class  
$$K_{\M_{0,14}}+\frac{1}{3}(B_2+B_5)+B_6+r_7B_7$$
where $r_7\in [1,4/3]$, is not of this form.

\vskip 10pt
\noindent {\tiny DEPARTMENT OF MATHEMATICS, UNIVERSITY OF MICHIGAN, ANN ARBOR, MI 48109-1109}
\vskip 3pt
{\footnotesize{\mbox{ } E-mail:{\ {\tt gfarkas@umich.edu}}
\vskip 5pt
\noindent {\tiny DEPARTMENT OF MATHEMATICS, UNIVERSITY OF MICHIGAN, ANN ARBOR, MI 48109-1109}
\vskip 3pt
{\footnotesize{\mbox{ } E-mail:{\ {\tt agibney@umich.edu}}


\begin{thebibliography}{[AC]}
       \bibitem[AC]{AC}
                E. Arbarello,\ M. Cornalba,\ 
                 \emph{Calculating cohomology groups of moduli spaces 
of curves via algebraic
                  geometry}, 
                 Inst. Hautes \'Etudes Sci. Publ. Math. No. 88(1998), 97-127.

      
        
         \bibitem[Fa1]{Fa1}
               C.  Faber,
                 \emph{Intersection-theoretical computations on $\M_g$},
Parameter Spaces (Warsaw 1994), 71-81, Banach Center Publ. 36, 1996.

         \bibitem[Fa2]{Fa2} 
               C. Faber, \emph{The nef cone of $\M_{0,6}$: a proof by inequalities only,} \ preprint.
        
         \bibitem[G]{g} A. Gibney, \emph{Fibrations of $\M_{g,n}$}, Ph.D. Thesis, University of Texas, 2000.
         \bibitem[GKM]{GKM}
               A.  Gibney,\ S. Keel,\ I. Morrison,
                 \emph{Towards the ample cone of $\M_{g,n}$}, to appear in J. Amer. Math. Soc.,
                math.AG/0006208.


\bibitem[HT]{ht} B. Hassett,\ Y. Tschinkel, \emph{On the effective cone of the moduli space of pointed rational curves}, math.AG/0110231.

         \bibitem[H]{H} B. Hunt,\ \emph{The geometry of some special arithmetic quotients}, Lecture Notes in Mathematics 1637, Springer 1996.

         \bibitem[HMo]{HMo}
               J.  Harris,\ I. Morrison, 
                 \emph{Moduli of curves},
                 Springer, 1998.
\bibitem[Kap]{Kap}
                M. Kapranov,
\emph{Veronese curves and the Grothendieck-Knudsen moduli space $\M_{0,n}$},
  J. of Algebraic Geometry, 2 (1993), 239-262.
         \bibitem[Ke]{Ke}
               S.  Keel,
                 \emph{Intersection theory on moduli spaces of $n$-pointed curves of genus zero}, Trans. Amer. Math. Soc. 330(1992), 545-574. 
         \bibitem[KMcK]{KMcK}
                S. Keel,\ J. McKernan,\ 
                 \emph{Contractible extremal rays on $\M_{0,n}$,}
                 math.AG/9607009.

\bibitem [K]{K}
           J.  Koll\'ar,\ 
\emph{Rational curves on algebraic varieties}, 
Springer 1996.
         

\bibitem [Ve]{Ve} P. Vermeire,\ \emph{A counterexample to Fulton's conjecture on $\M_{0,n}$}, preprint.
     

     
         \end{thebibliography}
\end{document}